\documentclass{elsarticle}

\usepackage{amssymb}
\usepackage{amsmath}

\usepackage[english,francais]{babel}

\newtheorem{theorem}{Theorem}[section]

\newtheorem{e-proposition}[theorem]{Proposition}

\newtheorem{e-definition}[theorem]{Definition\rm}

\def\og{\leavevmode\raise.3ex\hbox{$\scriptscriptstyle\langle\!\langle$~}}
\def\fg{\leavevmode\raise.3ex\hbox{~$\!\scriptscriptstyle\,\rangle\!\rangle$}}

\journal{the Acad\'emie des sciences}
\begin{document}
\begin{frontmatter}

\selectlanguage{english}
\title{Nonuniform Berry-Esseen bound for self-normalized martingales}


\selectlanguage{english}


\author[label1]{Songqi Wu}
\ead{sqwu@tju.edu.cn}
\author[label2]{Lingjie Kong}

\address[label1]{Center for Applied Mathematics,
Tianjin University, Tianjin, China}
\address[label2]{School of Mathematics, Tianjin University, Tianjin, China}



\begin{abstract}
\selectlanguage{english}
We give a nonuniform Berry-Esseen bound for self-normalized martingales, which bridges the gap between the results  of  Haeusler (1988) and Fan and Shao (2018). The bound coincides with the nonuniform Berry-Esseen bound of Haeusler and Joos (1988) for standardized martingales. As a consequence, a Berry-Esseen bound is obtained.
%

\end{abstract}
\end{frontmatter}

\section{Introduction}
Let $(X_i,\mathcal{F}_i)_{i=0,...,n}$ be a martingale differences defined on a probability space $(\Omega,\mathcal{F},\mathbf{P})$, where $X_0 = 0$ and $\{\emptyset,\Omega\} = \mathcal{F}_0 \subset ... \subset \mathcal{F}_n \subset \mathcal{F}$ are increasing $\sigma-$fields. Let $$S_0 = 0,\quad S_k = \sum_{i=1}^{k}X_i,\quad k=1,...,n.$$
Then $(S_k,\mathcal{F}_k)_{k=0,...,n}$ is a martingale. Without loss generality, assume that $\mathbf{E}S_n^2 = 1$, which means $S_n$ is standardized. Let $[S]$ and $\langle S \rangle$ be, respectively, the squared variance and the conditional variance of the martingale, that is $$[S]_0 = 0,\quad [S]_k = \sum_{i=1}^{k}X_i^2$$ and $$\langle S \rangle_0 = 0,\quad \langle S \rangle_k = \sum_{i=1}^{k}\mathbf{E}[X_i^2|\mathcal{F}_{i-1}],\quad k=1,...,n.$$ 

The absolute errors of normal approximations for martingales have been intensively studied; see, for instance, Heyde and Brown \cite{HB70}, Haeusler \cite{H88}, Haeusler and Joos \cite{HJ88}, El Machkouri and Ouchti \cite{M07}, Fan and Shao \cite{FS18}, Fan \cite{F19} and \cite{WMSF20}. Suppose that $\mathbf{E}|X_i|^{2p} \ \textless \ \infty$ for some $p \ \textgreater\ 1$ and all $i=1,...,n.$ Define
\begin{eqnarray}
N_n = \sum_{i=1}^{n}\mathbf{E}|X_i|^{2p} + \mathbf{E}|\langle S \rangle_n - 1|^p.
\end{eqnarray} Heyde and Brown \cite{HB70} proved the following Berry-Esseen bound 
\begin{eqnarray}\label{bou}
\sup_{x \in \mathbf{R}}\big|\mathbf{P}(S_n \le x) - \Phi(x)\big| \le C_pN_n^{1/(2p+1)},
\end{eqnarray}
where $p \in (1,2]$ and $C_p$ is a positive constant depending only on $p$. As a major advance in this direction, Haeusler \cite{H88} extended the result of Heyde and Brown \cite{HB70} from $p \in (1,2]$ to $p \in (1,\infty)$.
 Additionally, to justify that the bound (\ref{bou}) is asymptotically the best possible, he showed that there exists a martingale difference sequence $(\xi_k,\mathcal{F}_k)_{k=0,...,n}$, such that for $n$ large enough, $$\sup_{x \in \mathbf{R}}\big|\mathbf{P}(S_n \le x) - \Phi(x)\big|N_n^{-1/(2p+1)} \ge c_p,$$ where $c_p$ is a positive constant and does not depend on $n$. Later, Fan and Shao \cite{FS18} obtained that the following Berry-Esseen bounds for self-normalized and normalized martingales: for $p > 1$,
 \begin{eqnarray}\label{bouu}
 \sup_{x \in \mathbf{R}}\big|\mathbf{P}(S_n/\sqrt{[S]_n} \le x) - \Phi(x)\big| \le C_pN_n^{1/(2p+1)}
 \end{eqnarray}
 and
 \begin{eqnarray}\label{bouuu}
 \sup_{x \in \mathbf{R}}\big|\mathbf{P}(S_n/\sqrt{\langle S \rangle_n} \le x) - \Phi(x)\big| \le C_pN_n^{1/(2p+1)}.
 \end{eqnarray}
Notice that the last two bounds coincide with Haeusler's bound \cite{H88}. They also showed that the last two Berry-Esseen bounds are also the best possible in the spirit of Haeusler \cite{H88}. The first goal of this paper is to bridge the Berry-Esseen bounds of Haeusler \cite{H88}, Fan and Shao \cite{FS18}.
 
On the other hand and Haeusler and Joos \cite{HJ88} obtained the following nonuniform version of (\ref{bou}), that is for $p > 1$,
\begin{eqnarray}\label{nbou}
\big|\mathbf{P}(S_n \le x) - \Phi(x)\big| \le C_p\frac{N_n^{1/(2p+1)}}{1+|x|^{2p}}.
\end{eqnarray}
Clearly, the last bound improved the Berry-Esseen bounds by adding a factor decaying at a rate of polynomial. Inspired by the result of Haeusler and Joos \cite{HJ88}, 
we are interested in giving a nonuniform version of Berry-Esseen bound for Fan and Shao \cite{FS18}, which is the second goal of this paper.

All over the paper, $c$ and $C_p$, probably enabled with some indices, denote a common positive constant and a common positive constant depending only on $p$ respectively.

\section{Main results}
The following theorem bridges the Berry-Esseen bounds of Haeusler \cite{H88} and Fan and Shao \cite{FS18}.    
\begin{theorem}\label{th1}
	Assume that $\mathbf{E}|X_i|^{2p} \ \textless \ \infty$ for some $p \ \textgreater \ 1$ and all $i = 1,...,n$. Then there exists a constant $C_p$ denpending only on $p$ such that
	\begin{eqnarray}\label{th11}
	\sup_{x \in \mathbf{R}}\Big|\mathbf{P}\Big(\frac{S_n}{\sqrt{t_1\langle S \rangle_n + t_2[S]_n + (1 - t_1 - t_2)}} \le x\Big) - \Phi(x)\Big| \le C_pN_n^{1/(2p+1)},
	\end{eqnarray}
	where $t_1, t_2, (t_1 + t_2) \in [0,1].$
\end{theorem}

When $t_1 = 1$, the result (\ref{th11}) is exactly the Berry-Esseen bound for normalized martigale. When $t_2 = 1$, the result (\ref{th11}) reduces to the Berry-Esseen bound for self-normalized martingale. When $t_1 = t_2 = 0$, the result (\ref{th11}) becomes the Berry-Esseen bound of Haeusler \cite{H88}. Thus, the result bridges the Berry-Esseen bounds of Haeusler \cite{H88} and Fan and Shao \cite{FS18}.


Moreover, we also have the following nonuniform Berry-Esseen bound. 

\begin{theorem}\label{th2}
	Assume that $\mathbf{E}|X_i|^{2p} \ \textless \ \infty$ for some $p \ \textgreater \ 1$ and all $i = 1,...,n$. Then there exists a constant $C_p$ denpending only on $p$ such that
	\begin{eqnarray}\label{th21}
	\bigg|\mathbf{P}\bigg(\frac{S_n}{\sqrt{t_1\langle S \rangle_n + t_2[S]_n + (1 - t_1 - t_2)}}\le x\bigg) - \Phi(x)\bigg| \le C_{t_2,p}\frac{N_n^{1/(2p+1)}}{1+|x|^{2p}},
	\end{eqnarray}
	where $t_2 \in [0,1)$ and $t_1, (t_1 + t_2) \in [0,1].$
\end{theorem}

Notice that the bound (\ref{th21}) coincides with the nonuniform Berry-Esseen bound of Haeusler and Joos \cite{HJ88}.

\section{Proofs of theorems}

\subsection{Proof of Theorem \ref{th1}}\label{sec3}
When $N_n \ \textgreater \ 1$, the result is trivial, so we suppose $N_n \le 1$. First, we give an upper bound for $\mathbf{P}\Big(S_n \le x\sqrt{t_1\langle S \rangle_n + t_2[S]_n + t_3}\Big) - \Phi(x), x \le 0.$ Let $\epsilon_n \in (0,\frac{1}{2}]$ be a positive number and let $t_3 = 1 - t_1 - t_2$. It is easy to see that for $x \le 0$,
	\begin{eqnarray}\label{upbo}
	&\ &\mathbf{P}\Big(S_n \le x\sqrt{t_1\langle S \rangle_n + t_2[S]_n + t_3}\Big) - \Phi(x)\nonumber\\
	&\le&\mathbf{P}\Big(S_n \le x\sqrt{t_1(1 - \epsilon_n) + t_2[S]_n + t_3}\Big) - \Phi(x) + \mathbf{P}\Big(\langle S \rangle_n \le 1 - \epsilon_n\Big)\nonumber\\
		&\le&\mathbf{P}\Big(S_n \le x\sqrt{t_1(1 - \epsilon_n) + t_2(1 - \epsilon_n) + (1 - t_1 - t_2)}\Big)\nonumber\\
		&\ & -\ \Phi(x) + \mathbf{P}\Big(\langle S \rangle_n \le 1 - \epsilon_n\Big) + \mathbf{P}\Big([S]_n \le 1 - \epsilon_n \Big)\nonumber\\
		&=&\mathbf{P}\Big(S_n \le x\sqrt{1 - (t_1 + t_2)\epsilon_n}\Big)\nonumber\\
		&\ & -\ \Phi(x) + \mathbf{P}\Big(\langle S \rangle_n \le 1 - \epsilon_n\Big) + \mathbf{P}\Big([S]_n \le 1 - \epsilon_n\Big)\nonumber\\
		&=& J_1 + J_2 + J_3 + J_4,
		\end{eqnarray}
		where 
		%
		\begin{eqnarray}
		J_1 &=& \mathbf{P}\Big(S_n \le x\sqrt{1 - (t_1 + t_2)\epsilon_n}\Big) - \Phi\Big(x\sqrt{1 - (t_1 + t_2)\epsilon_n}\Big),\nonumber\\
		J_2 &=& \Phi\Big(x\sqrt{1 - (t_1 + t_2)\epsilon_n}\Big) - \Phi(x),\nonumber\\
		J_3 &=& \mathbf{P}\Big(\langle S \rangle_n \le 1 - \epsilon_n\Big),\nonumber\\
		J_4 &=& \mathbf{P}\Big([S]_n \le 1 - \epsilon_n\Big).\nonumber
		\end{eqnarray}
		%
		%
		Using Haeusler's inequality \cite{H88}, we obtain
		\begin{eqnarray}
		J_1 \le C_{p,1}\Big(\sum_{i=1}^{n}\mathbf{E}|X_i|^{2p} + \mathbf{E}|\langle S \rangle_n - 1|^p\Big)^{1/(2p+1)}.
		\end{eqnarray}
		Taking one-term Taylor's expansion and using the fact that $e^{-x^2/2}|x| \le 1$, we get
		\begin{eqnarray}
		J_2 &\le& c_1e^{-x^2/2}|x|\Big(1 - \sqrt{1 - (t_1 + t_2)\epsilon_n}\Big)\nonumber\\
		&\le& c_1e^{-x^2/2}|x|\Big(1 -  \sqrt{1 - (t_1 + t_2)\epsilon_n}\Big)\Big(1 + \sqrt{1 - (t_1 + t_2)\epsilon_n}\Big)\nonumber\\
		&\le&c_2\epsilon_n.
		\end{eqnarray}
		Using Markov's inequality \cite{H88}, we have
		\begin{eqnarray}
		J_3 &\le& \mathbf{P}\Big(|\langle S \rangle_n - 1 | \ge  \epsilon_n\Big)\nonumber\\&\le&\epsilon_n^{-p}\mathbf{E}|\langle S \rangle_n - 1|^p.
		\end{eqnarray}
		By the inequality (4.11) of Fan and Shao \cite{FS18}, we obtain for $p > 1$
		\begin{eqnarray}
		J_4 \le C_{p,2}\epsilon_n^{-p}\Big[\sum_{i=1}^{n}\mathbf{E}|X_i|^{2p} + \Big(\sum_{i=1}^{n}\mathbf{E}|X_i|^{2p}\Big)^{p/(2p-2)} + \mathbf{E}\big|\langle S \rangle_n - 1\big|^{p}\Big].
		\end{eqnarray}
		Returning to (\ref{upbo}), we have for $p \ \textgreater\ 1$,
		\begin{eqnarray}
		&\ &\mathbf{P}\Big(S_n \le x\sqrt{t_1\langle S \rangle_n + t_2[S]_n + t_3}\Big) - \Phi(x)\nonumber\\
		&\ &\quad \le  C_{p,1}\Big[\sum_{i=1}^{n}\mathbf{E}|X_i|^{2p} + \mathbf{E}\big|\langle S \rangle_n - 1\big|^{p}\Big]^{1/(2p+1)} + c_2\epsilon_n + \epsilon_n^{-p}\mathbf{E}\big|\langle S \rangle_n - 1\big|^{p}\nonumber\\
		&\ &\quad\quad +\ C_{p,2}\epsilon_n^{-p}\Big[\sum_{i=1}^{n}\mathbf{E}|X_i|^{2p} + \Big(\sum_{i=1}^{n}\mathbf{E}|X_i|^{2p}\Big)^{{p/(2p-2)}} + \mathbf{E}\big|\langle S \rangle_n - 1\big|^{p}\Big]\nonumber\\
		&\ &\quad \le  C_{p,1}\Big[\sum_{i=1}^{n}\mathbf{E}|X_i|^{2p} + \mathbf{E}\big|\langle S \rangle_n - 1\big|^{p}\Big]^{1/(2p+1)} + c_2\epsilon_n\nonumber\\
		&\ &\quad\quad +\ (C_{p,2} + 1)\epsilon_n^{-p}\Big[\sum_{i=1}^{n}\mathbf{E}|X_i|^{2p} + \Big(\sum_{i=1}^{n}\mathbf{E}|X_i|^{2p}\Big)^{{p/(2p-2)}}  + \mathbf{E}\big|\langle S \rangle_n - 1\big|^{p}\Big].\nonumber
		\end{eqnarray}
		Taking
		\begin{eqnarray}
		\epsilon_n = \Big[\sum_{i=1}^{n}\mathbf{E}|X_i|^{2p} + \Big(\sum_{i=1}^{n}\mathbf{E}|X_i|^{2p}\Big)^{{p/(2p-2)}} + \mathbf{E}\big|\langle S \rangle_n - 1\big|^{p}\Big]^{1/{(p+1)}}
		\end{eqnarray}
		and using the fact that $p/((2p-2)(p+1)) \ge 1/(2p+1)$, we observe for $x \le 0$ and $p\ \textgreater \ 1$,
		\begin{eqnarray}
		&\ &\mathbf{P}\Big(S_n \le x\sqrt{t_1\langle S \rangle_n + t_2[S]_n + t_3}\Big) - \Phi(x)\nonumber\\
		&\ &\quad \le  C_{p,1}\Big[\sum_{i=1}^{n}\mathbf{E}|X_i|^{2p} + \mathbf{E}\big|\langle S \rangle_n - 1\big|^{p}\Big]^{{1}/(2p+1)}\nonumber\\
		&\ &\quad\quad +\ C_{p,3}\Big[\sum_{i=1}^{n}\mathbf{E}|X_i|^{2p} + \Big(\sum_{i=1}^{n}\mathbf{E}|X_i|^{2p}\Big)^{{p/(2p-2)}} + \mathbf{E}\big|\langle S \rangle_n - 1\big|^{p}\Big]^{1/(p+1)}\nonumber\\
		&\ &\quad \le  C_{p,1}N_n^{{1}/(2p+1)} + C_{p,4}N_n^{1/(p+1)} + C_{p,4}\Big(\sum_{i=1}^{n}\mathbf{E}|X_i|^{2p}\Big)^{{p}/{((2p-2)(p+1))}}\nonumber\\
		&\ &\quad \le  C_{p,5}N_n^{{1}/(2p+1)} + C_{p,4}\Big(\sum_{i=1}^{n}\mathbf{E}|X_i|^{2p}\Big)^{{1}/(2p+1)}\nonumber\\
		&\ &\quad \le  C_{p,6}N_n^{{1}/(2p+1)}.
		\end{eqnarray}
		
		Next, we give a lower bound for $\mathbf{P}\Big(S_n \le x\sqrt{t_1\langle S \rangle_n + t_2[S]_n + t_3}\Big) - \Phi(x), x \le 0.$ Similarly, let $\epsilon_n \in (0,\frac{1}{2}]$ be a positive number and let $t_3 = 1 - t_1 - t_2$, 
		\begin{eqnarray}
		&\ &\mathbf{P}\Big(S_n \le x\sqrt{t_1\langle S \rangle_n + t_2[S]_n + t_3}\Big) - \Phi(x)\nonumber\\
		&\ge&\mathbf{P}\Big(S_n \le x\sqrt{t_1(1 + \epsilon_n) + t_2[S]_n + t_3}, \langle S \rangle_n \ \textless\ 1 + \epsilon_n\Big) - \Phi(x)\nonumber\\
		&\ge&\mathbf{P}\Big(S_n \le x\sqrt{t_1(1 + \epsilon_n) + t_2[S]_n + t_3}\Big) - \Phi(x) - \mathbf{P}\Big(\langle S \rangle_n \ge 1 + \epsilon_n\Big)\nonumber\\
		&\ge&\mathbf{P}\Big(S_n \le x\sqrt{t_1(1 + \epsilon_n) + t_2(1 + \epsilon_n)  + (1 - t_1 - t_2)}\Big)\nonumber\\
		&\ & -\ \Phi(x) - \mathbf{P}\Big(\langle S \rangle_n \ge 1 + \epsilon_n\Big)- \mathbf{P}\Big([S]_n \ge 1 + \epsilon_n\Big)\nonumber\\
			&=&\mathbf{P}\Big(S_n \le x\sqrt{(1 + (t_1 + t_2)\epsilon_n)}\Big)\nonumber\\
			&\ & -\ \Phi(x) - \mathbf{P}\Big(\langle S \rangle_n \ge 1 + \epsilon_n\Big) - \mathbf{P}\Big([S]_n \ge 1 + \epsilon_n\Big)\nonumber\\
			&=& J_5 + J_6 - J_7 - J_8,\nonumber
			\end{eqnarray}
			where
			\begin{eqnarray}
			J_5 &=& \mathbf{P}\Big(S_n \le x\sqrt{1 + (t_1 + t_2)\epsilon_n}\Big) - \Phi\Big( x\sqrt{1 + (t_1 + t_2)\epsilon_n}\Big),\nonumber\\
			J_6 &=& \Phi\Big( x\sqrt{1 + (t_1 + t_2)\epsilon_n}\Big) - \Phi(x),\nonumber\\
			J_7 &=& \mathbf{P}\Big(\langle S \rangle_n \ge 1 + \epsilon_n\Big),\nonumber\\
			J_8 &=& \mathbf{P}\Big([S]_n \ge 1 + \epsilon_n\Big).\nonumber
			\end{eqnarray}
			%
			Similarly, we can prove that for $x \le 0$ and $p\ \textgreater\ 1$,
			\begin{eqnarray}
			\mathbf{P}\Big(S_n \le x\sqrt{t_1\langle S \rangle_n + t_2[S]_n + t_3}\Big) - \Phi(x) \ge  -C_{p,10}N_n^\frac{1}{2p+1}.
			\end{eqnarray}
			Therefore, it follows that
			\begin{eqnarray}
			\sup_{x \le 0}\Big|\mathbf{P}\Big(S_n \le x\sqrt{t_1\langle S \rangle_n + t_2[S]_n + t_3}\Big) - \Phi(x)\Big| \le C_{p,9}N_n^\frac{1}{2p+1}.
			\end{eqnarray} 
			Note that $(-S_k,\mathcal{F}_k)_{k=0,...,n}$ is also a martingale. Applying the last inequality to $(-S_k,\mathcal{F}_k)_{k=0,...,n}$, we have
			\begin{eqnarray}\label{sim}
			&\ &\sup_{x \textgreater 0}\Big|\mathbf{P}\Big(S_n \le x\sqrt{t_1\langle S \rangle_n + t_2[S]_n + t_3}\Big) - \Phi(x)\Big|\nonumber\\&\ &\qquad =\sup_{x \textgreater 0}\Big|\mathbf{P}\Big(S_n \le x\sqrt{t_1\langle S \rangle_n + t_2[S]_n + t_3}\Big) - 1 + 1 - \Phi(x)\Big|\nonumber\\&\ &\qquad =\sup_{x \textgreater 0}\Big|\Phi(-x) - \mathbf{P}\Big(-S_n \le -x\sqrt{t_1\langle S \rangle_n + t_2[S]_n + t_3}\Big) \Big|\nonumber\\&\ &\qquad \le C_{p,11}N_n^\frac{1}{2p+1}.
			\end{eqnarray}
			Thus,
			\begin{eqnarray}
			\sup_{x \in \mathbf{R}}\Big|\mathbf{P}\Big(S_n \le x\sqrt{t_1\langle S \rangle_n + t_2[S]_n + t_3}\Big) - \Phi(x)\Big| \le C_{p,12}N_n^\frac{1}{2p+1}.
			\end{eqnarray}
\subsection{Proof of Theorem \ref{th2}}
When $N_n \ \textgreater \ 1$, the result is trivial, so we suppose $N_n \le 1$. First, we give a lower bound for $\mathbf{P}\Big(S_n \le x\sqrt{t_1\langle S \rangle_n + t_2[S]_n + t_3}\Big) -$$ \Phi(x), x \le 0.$ Let $\epsilon_n = \epsilon_n(x) \ge 0$ be a positive number and let $t_3 = 1 - t_1 - t_2$. It is easy to see that
\begin{eqnarray}
&\ &\mathbf{P}\Big(S_n \le x\sqrt{t_1\langle S \rangle_n + t_2[S]_n + t_3}\Big) - \Phi(x)\nonumber\\
&\ge&\mathbf{P}\Big(S_n \le x\sqrt{t_1(1 + \epsilon_n) + t_2[S]_n + t_3}\Big) - \Phi(x) - \mathbf{P}\Big(\langle S \rangle_n \ge 1 + \epsilon_n\Big)\nonumber\\
&\ge&\mathbf{P}\Big(S_n \le x\sqrt{t_1(1 + \epsilon_n) + t_2(1 + \epsilon_n) + (1 - t_1 - t_2)}\Big)\nonumber\\
&\ & -\ \Phi(x) - \mathbf{P}\Big(\langle S \rangle_n \ge 1 + \epsilon_n\Big) - \mathbf{P}\Big([S]_n \ge 1 + \epsilon_n\Big)\nonumber\\
&=&\mathbf{P}\Big(S_n \le x\sqrt{1 + (t_1 + t_2)\epsilon_n}\Big)\nonumber\\
&\ & -\ \Phi(x) - \mathbf{P}\Big(\langle S \rangle_n \ge 1 + \epsilon_n\Big) - \mathbf{P}\Big([S]_n \ge 1 + \epsilon_n\Big)\nonumber\\
&=& K_1(x) + K_2(x) - K_3(x) - K_4(x),\nonumber
\end{eqnarray}
where

\begin{eqnarray}
K_1(x) &=& \mathbf{P}\Big(S_n \le x\sqrt{1 + (t_1 + t_2)\epsilon_n}\Big) - \Phi\Big(x\sqrt{1 + (t_1 + t_2)\epsilon_n}\Big),\nonumber\\
K_2(x) &=& \Phi\Big(x\sqrt{1 + (t_1 + t_2)\epsilon_n}\Big) - \Phi(x),\nonumber\\
K_3(x) &=& \mathbf{P}\Big(\langle S \rangle_n \ge 1 + \epsilon_n\Big),\nonumber\\
K_4(x) &=& \mathbf{P}\Big([S]_n \ge 1 + \epsilon_n\Big).\nonumber
\end{eqnarray}
Using (\ref{nbou}), we get that
\begin{eqnarray}\label{lc1}
K_1(x) &\ge& -\frac{C_{p,1}N_n^{1/(2p+1)}}{1 + |x|^{2p}(1 + (t_1 + t_2)\epsilon_n)^p}\nonumber\\
&\ge& -\frac{C_{t_1,p,1}N_n^{1/(2p+1)}}{1 + |x|^{2p}}.
\end{eqnarray}
Taking one-term Taylor's expansion, we have
\begin{eqnarray}\label{lc2}
K_2(x) &\ge& -c_1e^{-x^2/2}|x|\Big(\sqrt{1 + (t_1 + t_2)\epsilon_n} - 1\Big)\nonumber\\
&\ge& -c_1e^{-x^2/2}|x|\Big(\sqrt{1 + (t_1 + t_2)\epsilon_n} - 1\Big)\Big(\sqrt{1 + (t_1 + t_2)\epsilon_n} + 1\Big)\nonumber\\
&\ge& -c_2e^{-x^2/2}|x|\epsilon_n.
\end{eqnarray}
By Markov's inequality, we obtain

\begin{eqnarray}\label{lc3}
K_3(x) &\le& \mathbf{P}\Big(\big|\langle S \rangle_n - 1\big| \ge \epsilon_n)\Big)\nonumber\\
&\le& \epsilon_n^{-p}\mathbf{E}\big|\langle S \rangle_n - 1\big|^p\nonumber\\
&\le& \epsilon_n^{-p}N_n\nonumber\\
&\le& \epsilon_n^{-p}N_n^{(p+1)/(2p+1)}.
\end{eqnarray}
Using the inequality (4.11) of Fan and Shao \cite{FS18} and the fact that $p/(2p-2) \ge (p+1)/(2p+1)$ and $N_n < 1$, we get for $p > 1$,
\begin{eqnarray}\label{lc4}
K_4(x) &\le& C_{p,2}\epsilon_n^{-p}\bigg( \sum_{i=1}^{n}\mathbf{E}|X_i|^{2p} +\Big(\sum_{i=1}^{n}\mathbf{E}|X_i|^{2p}\Big)^{p/(2p-2)} + \mathbf{E}\big|\langle S \rangle_n -1\big|^p\bigg)\nonumber\\&\le&C_{p,2}\epsilon_n^{-p}\Big(N_n + N_n^{p/(2p-2)}\Big)\nonumber\\
&\le&C_{p,2}\epsilon_n^{-p}\Big(N_n + N_n^{(p+1)/(2p+1)}\Big)\nonumber\\	&\le&2C_{p,2}\epsilon_n^{-p} N_n^{(p+1)/(2p+1)}.
\end{eqnarray} 

Then, combining (\ref{lc1}), (\ref{lc2}), (\ref{lc3}) and (\ref{lc4}) together, we have for all $x \le 0, p > 1$ and $\epsilon_n = \epsilon_n(x) \ge 0$,
\begin{eqnarray}
&\ &\mathbf{P}\Big(S_n \le x\sqrt{t_1\langle S \rangle_n + t_2[S]_n + t_3}\Big) - \Phi(x) \nonumber\\&\ &\ge -\frac{C_{t_1,p,1}N_n^{1/(2p+1)}}{1 + |x|^{2p}} - c_2e^{-x^2/2}|x|\epsilon_n -C_{p,2}\epsilon_n^{-p} N_n^{(p+1)/(2p+1)}.\nonumber
\end{eqnarray}
Taking
\begin{eqnarray}
\epsilon_n = \epsilon_n(x) = (1+x^2)N_n^{1/(2p+1)},
\end{eqnarray}
we have for all $x \le 0$,
\begin{eqnarray}\label{lb}
&\ &\mathbf{P}\Big(S_n \le x\sqrt{t_1\langle S \rangle_n + t_2[S]_n + t_3}\Big) - \Phi(x) \nonumber\\&\ &\ge -\frac{C_{t_1,p,1}N_n^{1/(2p+1)}}{1 + |x|^{2p}} - c_2e^{-x^2/2}|x|(1+x^2)N_n^{1/(2p+1)} -\frac{C_{p,2}N_n^{1/(2p+1)}}{(1+x^2)^p}\nonumber\\
&\ &\ge -\frac{C_{t_1,p,2}N_n^{1/(2p+1)}}{1 + |x|^{2p}} - \frac{c_2e^{-x^2/2}|x|(1+x^2)(1 + |x|^{2p})N_n^{1/(2p+1)}}{1 + |x|^{2p}}\nonumber\\
&\ &\ge  -\frac{C_{t_1,p,3}N_n^{1/(2p+1)}}{1 + |x|^{2p}} .
\end{eqnarray}

Next, we give an upper bound for $\mathbf{P}\Big(S_n \le x\sqrt{t_1\langle S \rangle_n + t_2[S]_n + t_3}\Big) - \Phi(x), x \le 0.$ Let $1 > \epsilon_n \ge 0$ and let $t_3 = 1 - t_1 - t_2$. We deduce that
\begin{eqnarray}
&\ &\mathbf{P}\Big(S_n \le x\sqrt{t_1\langle S \rangle_n + t_2[S]_n + t_3}\Big) - \Phi(x)\nonumber\\
&\le&\mathbf{P}\Big(S_n \le x\sqrt{t_1(1 - \epsilon_n) + t_2[S]_n + t_3}\Big) - \Phi(x) + \mathbf{P}\Big(\langle S \rangle_n \le 1 - \epsilon_n\Big)\nonumber\\
&\le&\mathbf{P}\Big(S_n \le x\sqrt{t_1(1 - \epsilon_n) + t_2(1 - \epsilon_n) + (1 - t_1 - t_2)}\Big)\nonumber\\
&\ & -\ \Phi(x) + \mathbf{P}\Big(\langle S \rangle_n \le 1 - \epsilon_n\Big) + \mathbf{P}\Big([S]_n \le 1 - \epsilon_n 
\Big)\nonumber\\
&=&\mathbf{P}\Big(S_n \le x\sqrt{1 - (t_1 + t_2)\epsilon_n}\Big)\nonumber\\
&\ & -\ \Phi(x) + \mathbf{P}\Big(\langle S \rangle_n \le 1 - \epsilon_n\Big) + \mathbf{P}\Big([S]_n \le 1 - \epsilon_n\Big)\nonumber\\
&=& K_5(x) + K_6(x) + K_7(x) + K_8(x),\nonumber
\end{eqnarray}
where
\begin{eqnarray}
K_5(x) &=& \mathbf{P}\Big(S_n \le x\sqrt{1 - (t_1 + t_2)\epsilon_n}\Big) - \Phi\Big(x\sqrt{1 - (t_1 + t_2)\epsilon_n}\Big),\nonumber\\
K_6(x) &=& \Phi\Big(x\sqrt{1 - (t_1 + t_2)\epsilon_n}\Big) - \Phi(x),\nonumber\\
K_7(x) &=& \mathbf{P}\Big(\langle S \rangle_n \le 1 - \epsilon_n\Big),\nonumber\\
K_8(x) &=& \mathbf{P}\Big([S]_n \le 1 - \epsilon_n\Big).\nonumber
\end{eqnarray}

Using an argument similar to the proof of (\ref{lb}), we know for all $x \le 0$ satisfying $\epsilon_n(x) \le 1$,
\begin{eqnarray}
\mathbf{P}\Big(S_n \le x\sqrt{t_1\langle S \rangle_n + t_2[S]_n + t_3}\Big) - \Phi(x) \le  \frac{C_{t_3,p,4}N_n^{1/(2p+1)}}{1 + |x|^{2p}} .
\end{eqnarray}
When $x \le 0$ satisfies $\epsilon_n(x) > 1$, that is $x < -\sqrt{N_n^{-1/(2p+1)} - 1}$, it is easy to see that

\begin{eqnarray}
&\ &\mathbf{P}\Big(S_n \le x\sqrt{t_1\langle S \rangle_n + t_2[S]_n + t_3}\Big) - \Phi(x) \nonumber\\&\le& \mathbf{P}\Big(S_n \le x\sqrt{t_3}\Big) - \Phi\big(x\sqrt{t_3}\big) + \Phi\big(x\sqrt{t_3}\big) - \Phi(x)\nonumber\\
&\le& \frac{C_pN_n^{1/(2p+1)}}{1+|x\sqrt{t_3}|^{2p}} + \frac{c_3e^{-x^2/2}|x|^3(1 - \sqrt{t_3})}{|x|^{2}}\nonumber\\
&\le &\frac{C_{t_3,p}N_n^{1/(2p+1)}}{1+|x|^{2p}} + \frac{c_3e^{-x^2/2}|x|^3(1 - \sqrt{t_3})}{N_n^{-1/(2p+1) }-1}\nonumber\\
&\le &\frac{C_{t_3,p}N_n^{1/(2p+1)}}{1+|x|^{2p}} + \frac{c_4e^{-x^2/2}|x|^3(1+|x|^{2p})N_n^{1/(2p+1)}}{1+|x|^{2p}}\nonumber\\
&\le &\frac{\hat{C}_{t_3,p}N_n^{1/(2p+1)}}{1+|x|^{2p}}. 
\end{eqnarray}	
Therefore, it follows that, for all $x \le 0$,
\begin{eqnarray}
\Big|\mathbf{P}\Big(S_n \le x\sqrt{t_1\langle S \rangle_n + t_2[S]_n + t_3}\Big) - \Phi(x)\Big| \le \frac{C_{t_3,p,5}N_n^{1/(2p+1)}}{1+|x|^{2p}}.
\end{eqnarray}
By an argument similar to that of (\ref{sim}), we get for all $x > 0$,
\begin{eqnarray}
&\ &\Big|\mathbf{P}\Big(S_n \le x\sqrt{t_1\langle S \rangle_n + t_2[S]_n + t_3}\Big) - \Phi(x)\Big|\nonumber\\&=&\Big|\mathbf{P}\Big(S_n \le x\sqrt{t_1\langle S \rangle_n + t_2[S]_n + t_3}\Big) - 1 + 1 - \Phi(x)\Big|\nonumber\\ &=&\Big|\Phi(-x) - \mathbf{P}\Big(-S_n \le -x\sqrt{t_1\langle S \rangle_n + t_2[S]_n + t_3}\Big) \Big|\nonumber\\ &\le& \frac{C_{t_3,p,5}N_n^{1/(2p+1)}}{1+|x|^{2p}}.
\end{eqnarray}
Thus,
\begin{eqnarray}
\Big|\mathbf{P}\Big(S_n \le x\sqrt{t_1\langle S \rangle_n + t_2[S]_n + t_3}\Big) - \Phi(x)\Big| \le \frac{C_{t_3,p,6}N_n^{1/(2p+1)}}{1+|x|^{2p}}.
\end{eqnarray}
This completes the proof of Theorem \ref{th2}

\selectlanguage{english}

\end{document}